\documentclass[runningheads]{llncs}
\usepackage[T1]{fontenc}
\usepackage{graphicx}
\usepackage{amsmath}
\usepackage{amssymb}

\usepackage{amsthm}
\usepackage{placeins}
\usepackage{hyperref}
\usepackage{floatrow}
\usepackage{color}

\theoremstyle{definition}

\begin{document}

\title{Accelerating Machine Learning Systems via Category Theory: Applications to Spherical Attention for Gene Regulatory Networks}

\titlerunning{Accelerating Machine Learning Systems via Category Theory}


\author{
Vincent Abbott\inst{1} 
\and 
Kotaro Kamiya\inst{2}
\and
Gerard Glowacki\inst{3} 
\and
Yu Atsumi\inst{2} 
\and
Gioele Zardini\inst{1}
\and
Yoshihiro Maruyama\inst{3} 
}

\authorrunning{V. Abbott et al.}

\institute{Laboratory for Information and Decision Systems, Massachusetts Institute of Technology, Cambridge, USA\\
\email{\{vtabbott,gzardini\}@mit.edu}
\and
SyntheticGestalt, Tokyo, Japan\\
\email{\{k.kamiya,y.atsumi\}@syntheticgestalt.com}
\and
School of Informatics, Nagoya University, Nagoya, Japan\\
\email{\{glowacki,maruyama\}@i.nagoya-u.ac.jp}
\thanks{Supported by JST (JPMJMS2033; JPMJPR24K9; JPMJFR206P).}
}

\maketitle

\begin{abstract}

How do we enable artificial intelligence models to improve themselves? This is central to exponentially improving generalized artificial intelligence models, which can improve their own architecture to handle new problem domains in an efficient manner that leverages the latest hardware. However, current automated compilation methods are poor, and efficient algorithms require years of human development. In this paper, we use neural circuit diagrams, based in category theory, to prove a general theorem related to deep learning algorithms, guide the development of a novel attention algorithm catered to the domain of gene regulatory networks, and produce a corresponding efficient kernel. The algorithm we propose, spherical attention, shows that neural circuit diagrams enable a principled and systematic method for reasoning about deep learning architectures and providing high-performance code. By replacing SoftMax with an $L^2$ norm as suggested by diagrams, it overcomes the special function unit bottleneck of standard attention while retaining the streaming property essential to high-performance. Our diagrammatically derived \textit{FlashSign} kernel achieves comparable performance to the state-of-the-art, fine-tuned FlashAttention algorithm on an A100, and $3.6\times$ the performance of PyTorch. Overall, this investigation shows neural circuit diagrams' suitability as a high-level framework for the automated development of efficient, novel artificial intelligence architectures.

\end{abstract}

\section{Introduction}
Deep learning architectures lack a systematic analytical framework. 
Traditional tools, such as linear algebra, fall short in capturing non-linearities, while graph-based representations and ad hoc diagrams fail to account for critical structural information such as broadcasting.
This omission blocks our understanding of resource usage and limits the development of efficient, parallelized implementations.
Consequently, automated compilation tools such as PyTorch's compile and Triton \cite{paszke_pytorch_2019} are typically constrained to basic, elementwise fusion.
More complex algorithms, including attention mechanisms~\cite{vaswani_attention_2017}, required years of engineering effort before GPU-optimized implementations became available~\cite{dao_flashattention_2022,dao_flashattention-2_2023,shah_flashattention-3_2024}.

\begin{figure}[tb]
    \centering
    \includegraphics[width=0.8\linewidth]{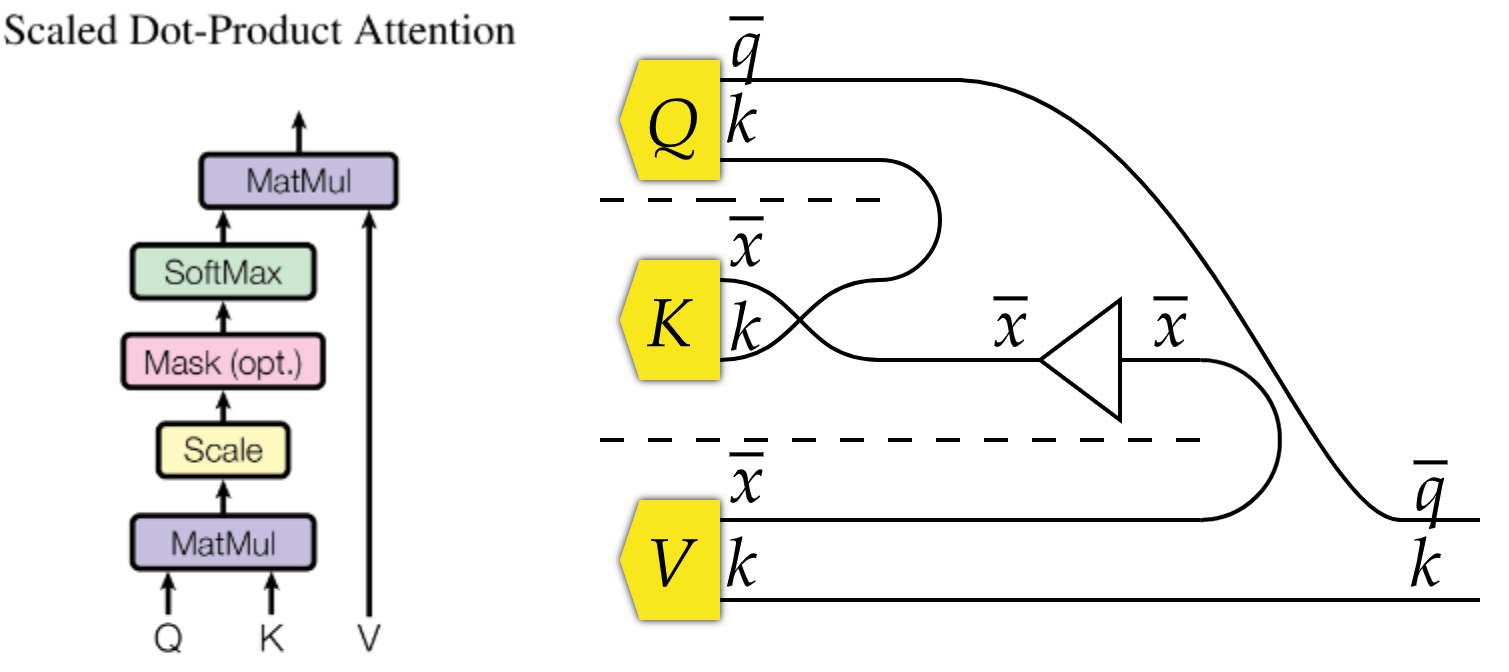}
    \caption{On the left is a traditional diagram from \textit{Attention is All You Need}\cite{vaswani_attention_2017}. On the right, we see a neural circuit diagram, which shows the details of axes and data sizes as well as broadcasting.}
    \label{fig:comparison}
\end{figure}

Neural circuit diagrams offer a compelling solution \cite{abbott_robust_2023,abbott_neural_2023,abbott_flashattention_2024}. 
They adapt monoidal string diagrams from category theory~\cite{selinger_survey_2009,xu_neural_2022}, encoding both data and operations of a deep learning model, with tensor axes represented as wires.
As illustrated in Figure \ref{fig:comparison}, dashed lines denote independence among operations or data structures, and broadcasting is expressed by enclosing operations with axis wires, which automatically determine input and output dimensions.
Linear contractions (e.g., dot products) appear as cups, while outer products are represented by open-ended dashed lines.

This formalism naturally accommodates non-linear operations, broadcasting, and linear algebraic structure.
More importantly, it exposes how computation is distributed across axes, enabling the derivation of efficient implementations aligned with the processor hierarchy.
These diagrams can be systematically mapped to low-level code, supporting precise performance modeling that includes memory hierarchies, cache-aware bandwidth estimates, and even overlapping clock-cycle analysis.

In unifying high-level algorithmic structure with hardware-aware execution detail, neural circuit diagrams make it possible to analyze, optimize, and redesign deep learning models within a single framework.
They provide the scaffolding needed for automating the discovery and refinement of learning algorithms, and as such, they lay the groundwork for building self-improving AI systems.

In order to show the utility of neural circuit diagrams, we will perform a principled investigation of a novel algorithm guided by the insights diagrammatic analysis. Specifically, the analysis of FlashAttention by neural circuit diagrams indicates that the exponential component of SoftMax is a bottleneck. The special function unit of GPUs has orders of magnitude fewer FLOPs available than matrix multiplication tensor cores or standard floating point cores \cite{dao_flashattention-2_2023}. For FP16 FlashAttention, the exponents take half as many clock cycles as the tensor core matrix multiplications. For FP8, they just take as many, completely bottlenecking the algorithm. Perfect overlapping is difficult, and there is reason to believe SoftMax has overhead beyond the special function unit operations \cite{abbott_flashattention_2024}, including warp-shuffling the maximum for numerical stability and unit conversion as the exponent can only occur in FP32. Diagrams, by revealing this information in clock cycle analysis, strongly hint at replacing the exponent with another operation.

However, the viability of attention as an algorithm requires the constituent steps, $\displaystyle QK^{T}$ matmul, the SoftMax normalization, and $\displaystyle SV$ matmul to be fused. This involves streaming input data in chunks and performing all steps without intermediate reads and writes to high-level memory. This prevents quadratic memory and bandwidth costs. The conditions under which fusion and streaming are possible is poorly understood using traditional methods, with fused, FlashAttention requiring five years to develop from the initial release of attention. Overcoming the exponential bottleneck requires the new algorithm to also be fused in the same manner. Helpfully, neural circuit diagrams allow the conditions for fusion to be systematically discovered and provide tools to prove that novel algorithms are fusable.

Furthermore, this work integrates the challenge of discovering new algorithms systematically, and in a potentially automated manner, with the field of categorical deep learning. Categorical deep learning \cite{gavranovic_categorical_2024} allows backpropagation \cite{fong_backprop_2019,cruttwell_categorical_2021}, \href{https://geometricdeeplearning.com/}{data-aware algorithms}, and architecture classes \cite{gavranovic_position_2024} to be formally described. More broadly, category theory has been applied to various fields of science, including physics \cite{baez_physics_2010}, logic \cite{abramsky_introduction_2010}, and others, revealing relationships between them and thus having promise for AGI \cite{abbott_category_2024}. It has also been applied to resource usage optimizations of complex engineered systems \cite{zardini_co-design_2023} and providing compilations \cite{wilson_categories_2022}, which means that a categorical formalism relates well to efficient implementation. This investigation shows that category theory, rather than being a pure abstract tool, can provide the structure for systematic derivation of deep learning algorithms.

In this paper, we follow the guidance of neural circuit diagrams to provide a novel attention algorithm that overcomes the exponential bottleneck, spherical attention, by replacing the SoftMax with an $\displaystyle L^{2}$-norm. We use neural circuit diagrams to prove a general theorem for streamable attention variants using alternative normalizers. We introduce spherical attention, and use its sign-preserving properties to provide an algorithm for gene regulatory networks. We provide a low-level implementation as guided by diagrams, achieving $3.6\times$ the performance of PyTorch and comparable throughput to FlashAttention on an A100. Overall, this shows the utility of neural circuit diagrams as a systematic and principled tool to analyze existing deep learning algorithms, reason about alternative architectures in a generalized manner, develop algorithms suited to novel domains, and provide efficient low-level implementations.

\section{Generalized Streamability}
Algorithms represented by neural circuit diagrams have the constituent axes of data and broadcasting clearly revealed. These parallelization details can be mapped onto execution in a processor hierarchy by recoloring wires to indicate presence at lower levels and relabelling to indicate the size and partition strategies they adopt. Group partitioning tiles an axis to be executed in parallel, and is feasible whenever a wire broadcasts the target operation. Streaming can be applied if an algorithm satisfies certain condition, and allows the final output to be calculated from chunks of the input data.

Streamability is accompanied by the \textit{fusion theorems}. These indicate that an algorithm remains streamable when composed with other operations and broadcasted over additional axes, as long as the streamed axis is preserved. The challenge of overcoming the exponential bottleneck of attention is developing an alternative which is streamable. Using the fusion theorems, we can prove that any normalization can replace the SoftMax operation to generate a streamable attention alternative.

We provide relevant definitions, a lemma, and a theorem below, which we prove in the Appendix \ref{apx:proofs}.

\begin{lemma}\label{lem:normalized_contraction}
    Normalized contractions are streamable.
\end{lemma}

\begin{figure}[h]
    \centering
\includegraphics[width=1\textwidth]{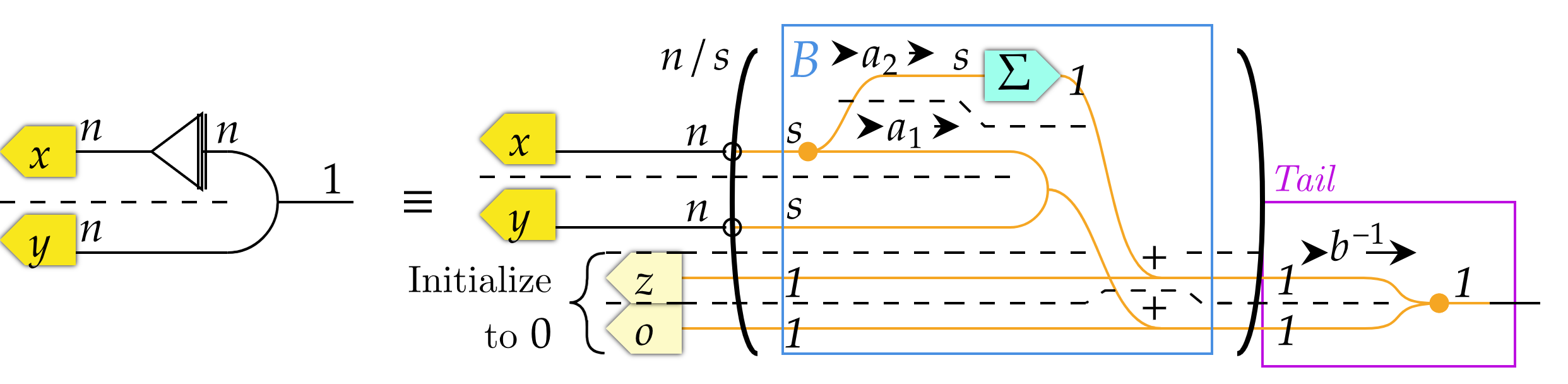}
    \caption{
    Generic normalized contraction, which we represent with a triple-lined triangle, can be expanded into a loop where the $n$-axis is partitioned into chunks of size $s$, and only $z$ and $o$ are maintained between chunks.
    }
    \label{fig:normalized_contraction}
\end{figure}
\FloatBarrier

\begin{theorem}\label{thm:alternative_attention}
Attention where the SoftMax is replaced by another normalization operation remains streamable, as shown in the diagram below;
\end{theorem}

\begin{figure}[h!]
    \floatbox[{\capbeside\thisfloatsetup{capbesideposition={left,top},capbesidewidth=0.5\textwidth}}]{figure}[\FBwidth]{
    \caption{
    An attention algorithm with a generic normalized contraction in place of SoftMax can have its $x$-axis streamed at a lower-level, resulting in the $s_x$ relabeling.
    }
    \label{fig:alternative_attention}}
    {\includegraphics[width=0.42\textwidth]{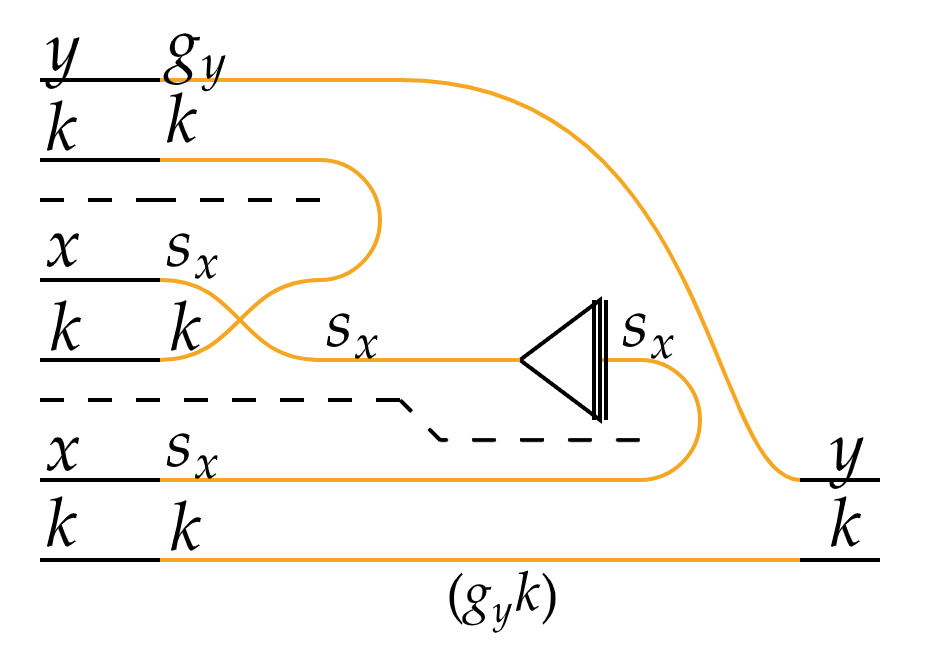}}
\end{figure}

\FloatBarrier

\section{Spherical Attention}
Theorem \ref{thm:alternative_attention} suggests we investigate forms of attention that overcome the exponential bottleneck by using an alternative normalization. Attention can be abstractly interpreted as a weighted linear combination of value vectors to attain a context-aware mixed vector output. This combination does not need a probabilistic interpretation, summing weights to one, like SoftMax.
There are various meaningful extensions of coefficients depending on the domain.
Here, we suggest spherical attention which uses an $\displaystyle L^{2}$-norm. Structurally, spherical attention has the advantage of being signed, allowing values to add or subtract from others, covering additional domains by providing a signalling rather than probabilistic interpretation of weights.
Furthermore, the exponential component of SoftMax bottlenecks throughput, especially at low quantization, which spherical attention overcomes by instead using fast FP16 cores.

As it matches the template of Figure \ref{fig:alternative_attention}, spherical attention can be efficiently implemented in a streamed manner. This reduces bandwidth requirements and removes the memory cost of storing $QK^T$, of size $y \times x$, which is quadratic in the standard case of self-attention. In contrast, a streamed attention algorithm only needs to store $g_y \times s_x$ data on each low-level chip (SM). These are configurable partition sizes, avoiding a memory limit with input size.

\subsection{Application to Gene Regulatory Networks}

Regulatory networks are graphs that consist of vertices and edges which indicate up- or downregulation . Gene regulatory networks (GRNs) \cite{aduddell_compositional_2024} are a concrete example, wherein genes up- or downregulate others to exhibit complex behaviour. A machine learning model can assist in estimating an underlying model from gene expression data, however, this domain requires further architectural features. Signed attention, by allowing positive and negative signals between tokens, provides the core of the structure.

Furthermore, genes can appear multiple times, amplifying their effect, meaning they have a bag structure. We set the token indexes to each refer to a specific gene, and multiplicity information for each gene is provided as a model input. This information is weaved into the attention model by scalar multiplication of keys at each layer. Input data is the multiplicity of each gene, capturing a snapshot of a single cell's state, and outputs are some concrete property, such as the cell identity. Layers implement successive regulatory networks.

The query-key attention scores are interpreted as the key gene's regulation of the query gene, which dictates how much of its value vector is to be added. Positive sign indicates up-regulation, while a negative sign indicates down-regulation. As the model is learned, genes' up- or downregulation of others is captured in $\displaystyle QK^{T}$ adjacency matrices.

The overall model design is shown below, in Figure \ref{fig:full_model}. The architectural features are explained in Appendix \ref{apx:full_model}.

\begin{figure}[h]
    \centering
\includegraphics[width=1\textwidth]{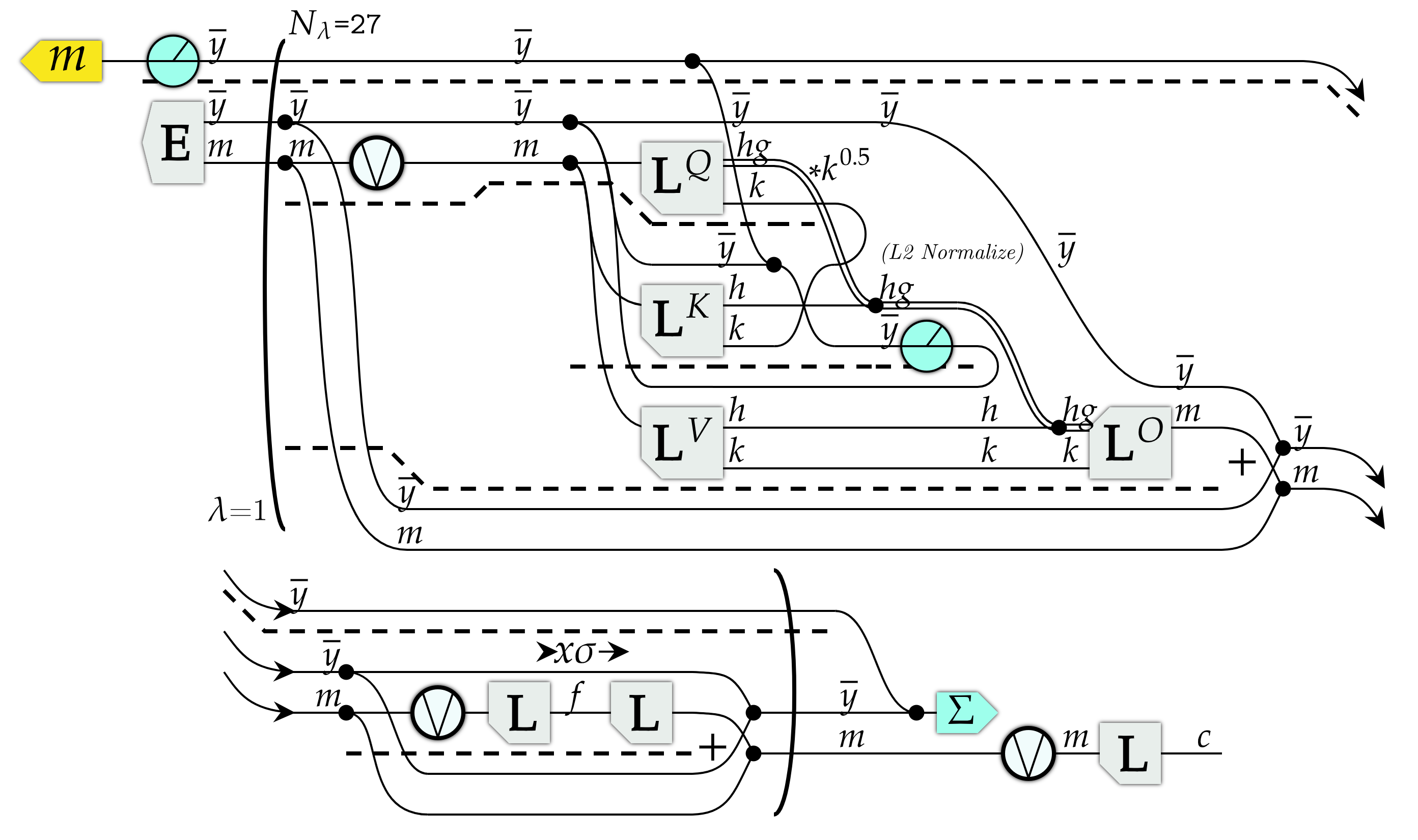}
    \caption{
    The full model, where spherical attention is combined with iterative feed-forward layers, can be diagrammatically expressed.
    }
    \label{fig:full_model}
\end{figure}
\FloatBarrier

\section{Low-Level Algorithm}\label{sec:lowlevel}
\begin{figure}[h!]
    \centering
    \includegraphics[width=1\linewidth]{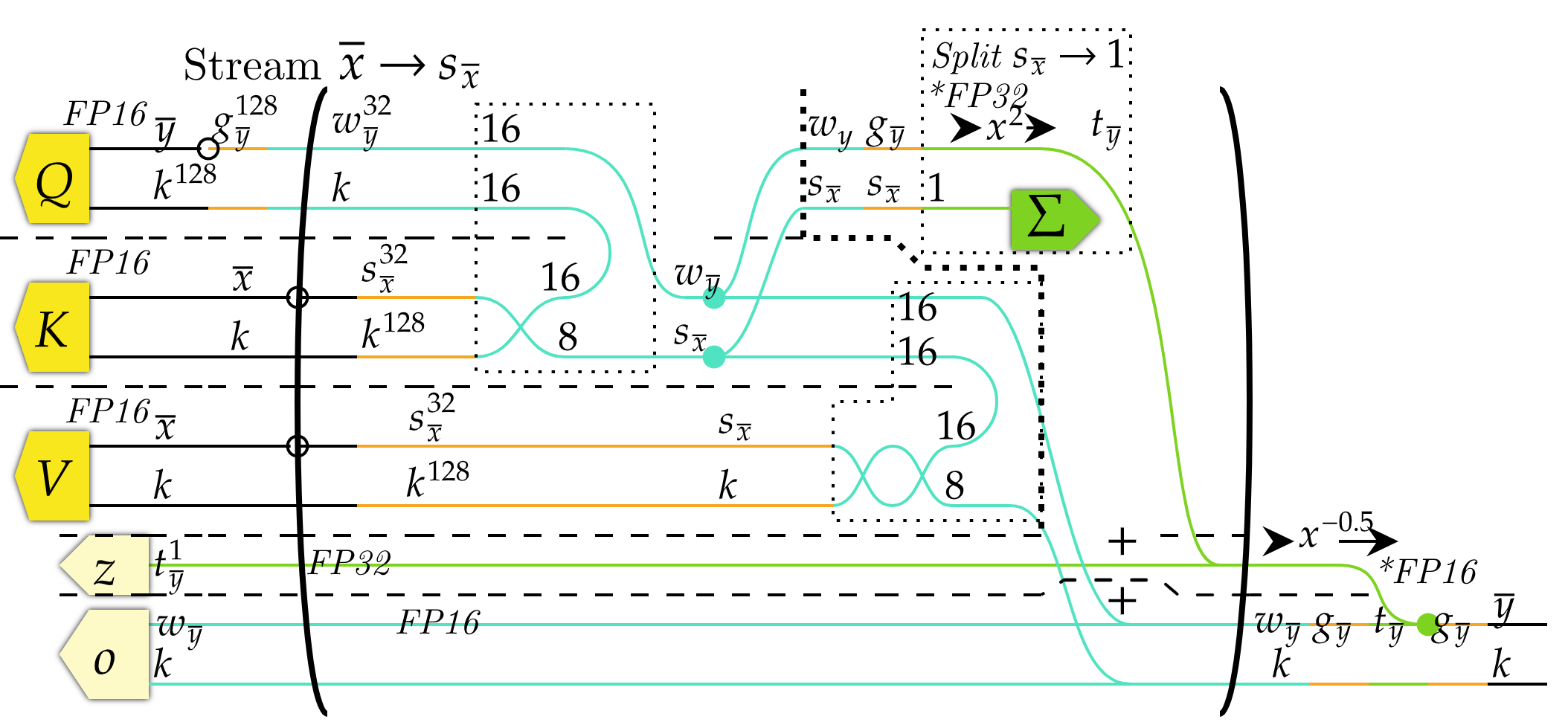}
    \caption{Low-level diagram of spherical attention, adapting Figure \ref{fig:attention_expansion} to derive diagrammatic pseudocode.}
    \label{fig:lowlevel}
\end{figure}

In Figure \ref{fig:lowlevel}, we show our low-level diagrammatic pseudocode. Variables are arrays of numbers at some quantization. They are represented on diagrams as a series of wires, labelled with each axis, with the lowest stride (the major axis) at the bottom. Wires are labelled with axis sizes, or tags indicating how an axis is split at various levels. Distinct variables are separated by dashed lines. We can tag these segments with the quantization used, and change quantization by adding an asterisk. Operations are denoted by pictograms, with their input variable size to the left and output to the right. Broadcasting involves drawing wires over operations, which naturally updates the size of input/output variables. Elementwise operations have input and output sizes of $\displaystyle 1$. Axes of size $\displaystyle 1$ can be discarded, meaning elementwise operations appear as floating arrows. These rules encompass the diagrams shown in previous figures. By revealing this parallelization information, deriving a low-level kernel is possible.

We assign colors to operations depending on the level at which they operate. We need to move data to on-chip SMEM, indicated orange. This allows us to move data to the registers of threads (green), or fragment data for tensor core computations (teal). Both use the same register RMEM memory. However, tensor core operations require fragmented memory which makes general-purpose operations difficult. As a result, we treat blue as a distinct memory space, requiring transfer through orange to get to green. Data is recursively tiled at each level. GMEM (black wires) have full axes, which are split into blocks at the SMEM level indicated by $\displaystyle g_{\square }$, split into warps for tensor cores $\displaystyle w_{\square }$, and split between threads for registers $\displaystyle t_{\square }$. When data moves between levels with different tile sizes, this represents partitioning or concatenating that axis. The main loop of the algorithm is encompassed within the parantheses, wherein $\displaystyle s_{\overline{x}}$ is split.

A dotted box indicates a linear operation which is split into a loop. Within a box, axes are relabelled. This is the per-iteration size. An input axis can change colors, indicating that it is loaded in chunks to lower levels, reducing memory usage. Tensor cores operate on tiles of the input data, so are typically split in such a manner. Tensor cores effectively implement $\displaystyle AB^{T}$, taking a $\displaystyle M\times K$ and $\displaystyle N\times K$ matrix to produce a $\displaystyle M\times N$ matrix by contracting the $\displaystyle K$ dimension. Therefore, to implement $\displaystyle SV$ we must transpose $\displaystyle V$ when loading into the tensor cores. This is indicated by the double transpose; the first transpose is intentionally implemented while loading $\displaystyle V$ into tensor cores, and the second is implicit in how tensor cores operate and are diagrammed. Overall, this results in a double transpose, functionally equivalent to an identity.

We can apply these rules to the diagram from Figure \ref{fig:attention_expansion}, to derive the low-level pseudocode diagram in Figure \ref{fig:lowlevel}. Following the procedure above, we colored tensor core data and operations blue, and general-purpose thread operations green. Data is colored orange as it transfers through SMEM. Tile sizes for levels have been assigned by labelling. Linear operations, including matrix multiplication and summation, have been split via dotted boxes, reducing memory usage. We have also added a thick dotted line to indicate a necessary synchronisation. Axes sizes are indicated by superscripts.

\FloatBarrier

\subsection{Results}

We implemented our algorithm in CUDA C++ using inline PTX, giving us fine-grained control over asynchronous GMEM to SMEM writes and tensor cores needed to pad out the $K$ and $V$ buffers in SMEM. These were integrated into PyTorch using inline loading of the CUDA source code.

For our testing and benchmarks, $Q$, $K$, and $V$ had standard normal distributions. We chose multiples of $13824=108\times128$ as our sequence lengths. This allows our algorithm to be efficiently wave-quantized on an A100. Head dimension $128$ and FP16 quantization matches the algorithm we expressed in Section \ref{sec:lowlevel}.

\begin{figure}
    \centering
    \includegraphics[width=1\linewidth]{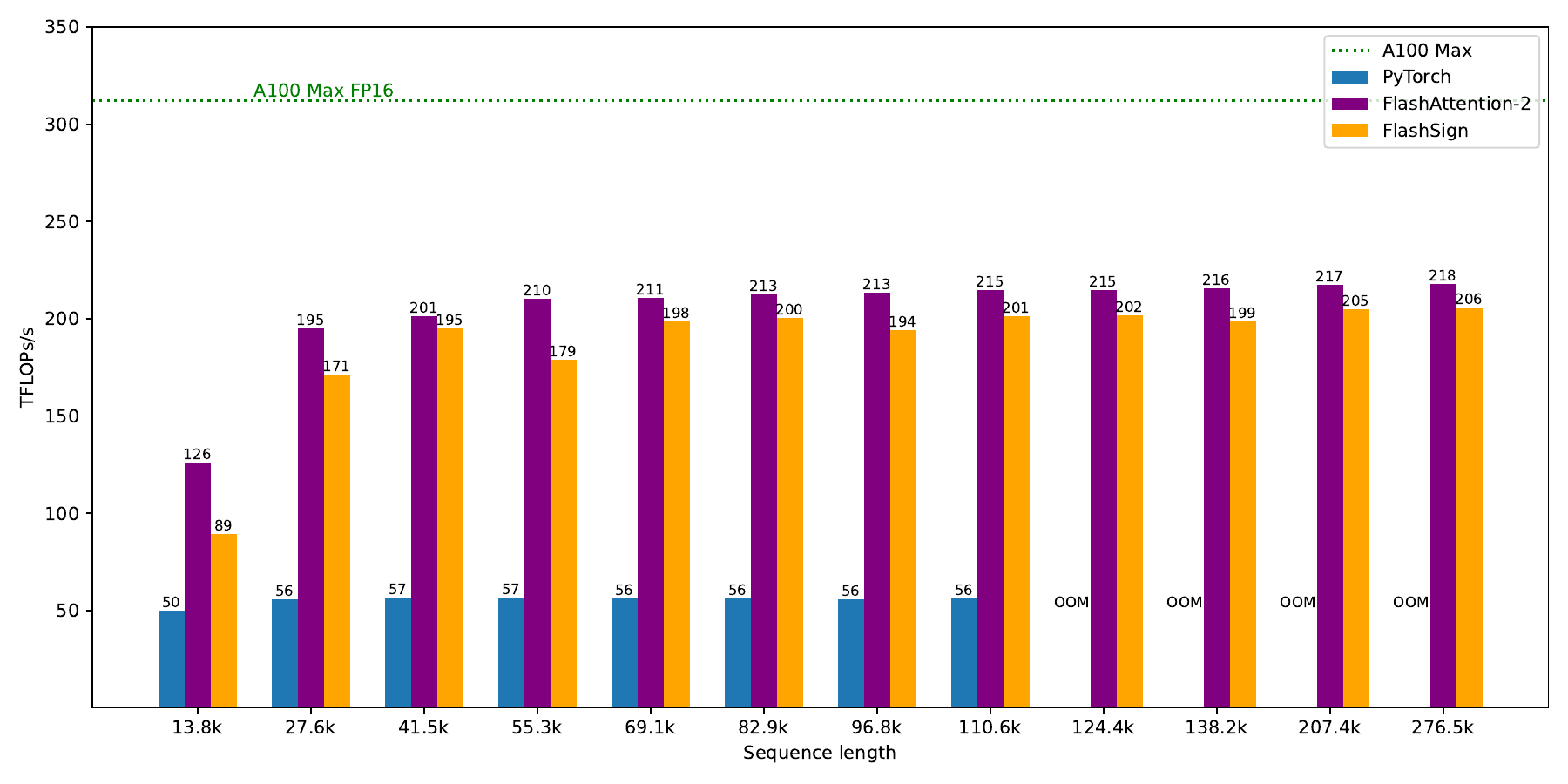}
    \caption{Results comparing our algorithm to state-of-the-art FlashAttention}
    \label{fig:benchmark}
\end{figure}
\FloatBarrier

For large sequence lengths, we achieve around 200 TFLOP/s, or 64\% utilization of an A100's 312 TFLOP/s of FP16 tensor core operations. This equates to a speedup of up to 3.60$\times$ compared to PyTorch. We also see that PyTorch runs into out-of-memory (OOM) errors at large sequence lengths.

For comparison, we also ran \href{https://github.com/Dao-AILab/flash-attention}{FlashAttention-2} \cite{dao_flashattention-2_2023} on our system at the same sequence lengths. FlashAttention-2 is a state-of-the-art algorithm, with substantial hyperparameter fine-tuning and advanced features, such as warp-shuffling to avoid intermediate transfers to SMEM. Diagrammatically, this is equivalent to moving directly between green and blue levels. Despite this, our algorithm achieves within 5\% of their performance at a sequence length of 41.5K.

\section{Conclusions} 
In this work, we demonstrated how diagrammatic representations can be used both to derive structural properties of deep learning architectures and to guide the construction of efficient low-level kernels.
This stands in contrast to conventional approaches, which struggle to bridge the gap between high-level model design and hardware-aware implementation.
Our findings suggest that the representational and analytical capabilities of neural circuit diagrams offer a powerful foundation for the systematic improvement, and potentially the automated discovery, of deep learning algorithms, a critical step toward achieving AGI.

Looking ahead, we aim to automate the diagrammatic framework to enable the procedural generation of learning algorithms.
This would allow neural circuit diagrams to support rapid analysis, optimization, and discovery at scale.
Achieving this requires a formal data structure capable of syntactically encoding diagrammatic constructs, currently relying on human visual interpretation.
The categorical foundations of neural circuit diagrams make such a representation, and its corresponding compilation pipeline, both natural and feasible.
Ultimately, automating the tools of neural circuit diagrams can form the basis of self-discovering and improving algorithms.

\bibliographystyle{plain}
\bibliography{references}

\begin{thebibliography}{10}

\bibitem{abbott_neural_2023}
Vincent Abbott.
\newblock Neural {Circuit} {Diagrams}: {Robust} {Diagrams} for the {Communication}, {Implementation}, and {Analysis} of {Deep} {Learning} {Architectures}.
\newblock {\em Accepted to Transactions on Machine Learning Research}, July 2023.

\bibitem{abbott_robust_2023}
Vincent Abbott.
\newblock {\em Robust {Diagrams} for {Deep} {Learning} {Architectures}: {Applications} and {Theory}}.
\newblock Honours {Thesis}, The Australian National University, Canberra, October 2023.

\bibitem{abbott_category_2024}
Vincent Abbott, Tom Xu, and Yoshihiro Maruyama.
\newblock Category {Theory} for {Artificial} {General} {Intelligence}.
\newblock In {\em Artificial {General} {Intelligence}: 17th {International} {Conference}, {AGI} 2024, {Seattle}, {WA}, {USA}, {August} 13–16, 2024, {Proceedings}}, pages 119--129, Berlin, Heidelberg, August 2024. Springer-Verlag.

\bibitem{abbott_flashattention_2024}
Vincent Abbott and Gioele Zardini.
\newblock {FlashAttention} on a {Napkin}: {A} {Diagrammatic} {Approach} to {Deep} {Learning} {IO}-{Awareness}.
\newblock {\em Transactions on Machine Learning Research}, December 2024.

\bibitem{abramsky_introduction_2010}
Samson Abramsky and Nikos Tzevelekos.
\newblock Introduction to {Categories} and {Categorical} {Logic}.
\newblock volume 813, pages 3--94. 2010.
\newblock arXiv:1102.1313 [cs, math].

\bibitem{aduddell_compositional_2024}
Rebekah Aduddell, James Fairbanks, Amit Kumar, Pablo~S. Ocal, Evan Patterson, and Brandon~T. Shapiro.
\newblock A compositional account of motifs, mechanisms, and dynamics in biochemical regulatory networks.
\newblock {\em Compositionality}, Volume 6 (2024), May 2024.
\newblock Publisher: Episciences.org.

\bibitem{ainslie_gqa_2023}
Joshua Ainslie, James Lee-Thorp, Michiel~de Jong, Yury Zemlyanskiy, Federico Lebrón, and Sumit Sanghai.
\newblock {GQA}: {Training} {Generalized} {Multi}-{Query} {Transformer} {Models} from {Multi}-{Head} {Checkpoints}, December 2023.
\newblock arXiv:2305.13245.

\bibitem{baez_physics_2010}
John~C. Baez and Mike Stay.
\newblock Physics, {Topology}, {Logic} and {Computation}: {A} {Rosetta} {Stone}.
\newblock volume 813, pages 95--172. 2010.
\newblock arXiv:0903.0340 [quant-ph].

\bibitem{cruttwell_categorical_2021}
G.~S.~H. Cruttwell, Bruno Gavranović, Neil Ghani, Paul Wilson, and Fabio Zanasi.
\newblock Categorical {Foundations} of {Gradient}-{Based} {Learning}, July 2021.
\newblock arXiv:2103.01931 [cs, math].

\bibitem{dao_flashattention-2_2023}
Tri Dao.
\newblock {FlashAttention}-2: {Faster} {Attention} with {Better} {Parallelism} and {Work} {Partitioning}, July 2023.
\newblock arXiv:2307.08691.

\bibitem{dao_flashattention_2022}
Tri Dao, Daniel~Y. Fu, Stefano Ermon, Atri Rudra, and Christopher Ré.
\newblock {FlashAttention}: {Fast} and {Memory}-{Efficient} {Exact} {Attention} with {IO}-{Awareness}, June 2022.
\newblock arXiv:2205.14135 [cs].

\bibitem{deepseek-ai_deepseek-v3_2025}
DeepSeek-AI.
\newblock {DeepSeek}-{V3} {Technical} {Report}, February 2025.
\newblock arXiv:2412.19437 [cs].

\bibitem{fong_backprop_2019}
Brendan Fong, David~I. Spivak, and Rémy Tuyéras.
\newblock Backprop as {Functor}: {A} compositional perspective on supervised learning.
\newblock In {\em 34th {Annual} {ACM}/{IEEE} {Symposium} on {Logic} in {Computer} {Science}, {LICS} 2019, {Vancouver}, {BC}, {Canada}, {June} 24-27, 2019}, pages 1--13. IEEE, 2019.

\bibitem{gavranovic_categorical_2024}
Bruno Gavranović, Paul Lessard, Andrew Dudzik, Tamara von Glehn, João G.~M. Araújo, and Petar Veličković.
\newblock Categorical {Deep} {Learning}: {An} {Algebraic} {Theory} of {Architectures}, February 2024.
\newblock arXiv:2402.15332 [cs, math, stat].

\bibitem{gavranovic_position_2024}
Bruno Gavranović, Paul Lessard, Andrew Dudzik, Tamara von Glehn, João G.~M. Araújo, and Petar Veličković.
\newblock Position: {Categorical} {Deep} {Learning} is an {Algebraic} {Theory} of {All} {Architectures}, February 2024.

\bibitem{he_deep_2015}
Kaiming He, Xiangyu Zhang, Shaoqing Ren, and Jian Sun.
\newblock Deep {Residual} {Learning} for {Image} {Recognition}.
\newblock {\em CoRR}, abs/1512.03385, 2015.
\newblock arXiv: 1512.03385.

\bibitem{he_identity_2016}
Kaiming He, Xiangyu Zhang, Shaoqing Ren, and Jian Sun.
\newblock Identity {Mappings} in {Deep} {Residual} {Networks}.
\newblock In Bastian Leibe, Jiri Matas, Nicu Sebe, and Max Welling, editors, {\em Computer {Vision} - {ECCV} 2016 - 14th {European} {Conference}, {Amsterdam}, {The} {Netherlands}, {October} 11-14, 2016, {Proceedings}, {Part} {IV}}, volume 9908 of {\em Lecture {Notes} in {Computer} {Science}}, pages 630--645. Springer, 2016.

\bibitem{jiang_mixtral_2024}
Albert~Q. Jiang, Alexandre Sablayrolles, Antoine Roux, Arthur Mensch, Blanche Savary, Chris Bamford, Devendra~Singh Chaplot, Diego de~las Casas, Emma~Bou Hanna, Florian Bressand, Gianna Lengyel, Guillaume Bour, Guillaume Lample, Lélio~Renard Lavaud, Lucile Saulnier, Marie-Anne Lachaux, Pierre Stock, Sandeep Subramanian, Sophia Yang, Szymon Antoniak, Teven~Le Scao, Théophile Gervet, Thibaut Lavril, Thomas Wang, Timothée Lacroix, and William~El Sayed.
\newblock Mixtral of {Experts}, January 2024.
\newblock arXiv:2401.04088.

\bibitem{paszke_pytorch_2019}
Adam Paszke, Sam Gross, Francisco Massa, Adam Lerer, James Bradbury, Gregory Chanan, Trevor Killeen, Zeming Lin, Natalia Gimelshein, Luca Antiga, Alban Desmaison, Andreas Köpf, Edward Yang, Zach DeVito, Martin Raison, Alykhan Tejani, Sasank Chilamkurthy, Benoit Steiner, Lu~Fang, Junjie Bai, and Soumith Chintala.
\newblock {PyTorch}: {An} {Imperative} {Style}, {High}-{Performance} {Deep} {Learning} {Library}, December 2019.
\newblock arXiv:1912.01703.

\bibitem{selinger_survey_2009}
Peter Selinger.
\newblock A survey of graphical languages for monoidal categories, August 2009.

\bibitem{shah_flashattention-3_2024}
Jay Shah, Ganesh Bikshandi, Ying Zhang, Vijay Thakkar, Pradeep Ramani, and Tri Dao.
\newblock {FlashAttention}-3: {Fast} and {Accurate} {Attention} with {Asynchrony} and {Low}-precision, July 2024.
\newblock arXiv:2407.08608.

\bibitem{vaswani_attention_2017}
Ashish Vaswani, Noam Shazeer, Niki Parmar, Jakob Uszkoreit, Llion Jones, Aidan~N. Gomez, Lukasz Kaiser, and Illia Polosukhin.
\newblock Attention is {All} you {Need}.
\newblock In Isabelle Guyon, Ulrike~von Luxburg, Samy Bengio, Hanna~M. Wallach, Rob Fergus, S.~V.~N. Vishwanathan, and Roman Garnett, editors, {\em Advances in {Neural} {Information} {Processing} {Systems} 30: {Annual} {Conference} on {Neural} {Information} {Processing} {Systems} 2017, {December} 4-9, 2017, {Long} {Beach}, {CA}, {USA}}, pages 5998--6008, 2017.

\bibitem{wilson_categories_2022}
Paul Wilson and Fabio Zanasi.
\newblock Categories of {Differentiable} {Polynomial} {Circuits} for {Machine} {Learning}, May 2022.
\newblock arXiv:2203.06430 [cs, math].

\bibitem{xu_neural_2022}
Tom Xu and Yoshihiro Maruyama.
\newblock Neural {String} {Diagrams}: {A} {Universal} {Modelling} {Language} for {Categorical} {Deep} {Learning}.
\newblock In Ben Goertzel, Matthew Iklé, and Alexey Potapov, editors, {\em Artificial {General} {Intelligence}}, Lecture {Notes} in {Computer} {Science}, pages 306--315, Cham, 2022. Springer International Publishing.

\bibitem{zardini_co-design_2023}
Gioele Zardini.
\newblock {\em Co-{Design} of {Complex} {Systems}: {From} {Autonomy} to {Future} {Mobility} {Systems}}.
\newblock Doctoral {Thesis}, ETH Zurich, 2023.
\newblock Accepted: 2023-12-19T10:03:57Z.

\end{thebibliography}
\appendix
\section{Appendix}
\subsection{Generalized Streamability}\label{apx:proofs}
Streamability is a property of functions that allows them to be computed with minimal memory on lower-level chips by streaming data in blocks and accumulating with partially computed outputs \cite{abbott_flashattention_2024}. Mathematically, this can be represented by disaggregating concatenated vectors, which maps to reduced memory usage required to compute the function.

\begin{definition}[Streamable Function]
    A polymorphic function $\displaystyle f:X^{n}\rightarrow Y$ defined for any $\displaystyle n\in \mathbb{N}$ is \textit{\textbf{streamable} }if there exists a polymorphic accumulator $\displaystyle B:Y\times X^{n}\rightarrow Y$ such that $\displaystyle B( f(\mathbf{x}) ,\mathbf{y}) =f(\mathbf{x} \oplus \mathbf{y})$, where $\displaystyle \oplus $ is concatenation along the $\displaystyle n$-axis. If $\displaystyle f:X^{n}\rightarrow Y$ is streamable and we have $\displaystyle t:Y\rightarrow Z$, then $\displaystyle t\circ f:X^{n}\rightarrow Z$ is also considered streamable. Streamable functions require little memory at lower levels to compute, as only a limited $\displaystyle X^{n}$ and $\displaystyle Y$ are required to be on lower levels.
\end{definition}

\begin{definition}[Normalized Contraction]
 A \textit{\textbf{normalized operation }}$\displaystyle N:\mathbb{R}^{n}\rightarrow \mathbb{R}^{n}$ uses a pair of activation functions $\displaystyle a_{1}$, $\displaystyle a_{2}$ and an aggregator function $\displaystyle b$ to provide $\displaystyle N(\mathbf{x})_{i} =a_{1}( x_{i}) /b( \Sigma _{j} a_{2}( x_{j}))$. A \textbf{\textit{normalized contraction}} $\displaystyle NC:\mathbb{R}^{n} \times \mathbb{R}^{n} \simeq \left(\mathbb{R}^{2}\right)^{n}\rightarrow \mathbb{R}$ follows normalization with a linear contaction, $\displaystyle NC(\mathbf{x} ,\mathbf{y}) =N(\mathbf{x}) \cdot \mathbf{y} =( \Sigma _{i} a_{1}( x_{i}) y_{i}) /b( \Sigma _{i} a_{2}( x_{i}))$.
\end{definition}

\begin{lemma}
    Normalized contractions are streamable.
\end{lemma}

\textit{Proof}. Set up $\displaystyle B:\mathbb{R}^{2} \times \left(\mathbb{R}^{n}\right)^{2}\rightarrow \mathbb{R}^{2}$ so that $\displaystyle B(( o,z) ,(\mathbf{x} ,\mathbf{y})) \mapsto ( o+\Sigma _{i} a_{1}( x_{i}) \cdot y_{i} ,z+\Sigma _{i} a_{2}( x_{i}))$. Then, we see that;
\begin{align*}
B( B(( o,z) ,(\mathbf{x} ,\mathbf{y})) ,\ (\mathbf{x} ',\mathbf{y} ')) & =B(( o+\Sigma _{i} a_{1}( x_{i}) \cdot y_{i} ,z+\Sigma _{i} a_{2}( x_{i})) ,\ (\mathbf{x} ',\mathbf{y} '))\\
 & =( o+\Sigma _{i} a_{1}( x_{i}) \cdot y_{i} +\Sigma _{j} a_{2}( x_{j}) \cdot y_{j} ,z+\Sigma _{i} a_{2}( x_{i}) +\Sigma _{j} a_{2}( x_{j}))\\
 & =B(( o,z) ,(\mathbf{x} \oplus \mathbf{x} ',\mathbf{y} \oplus \mathbf{y} '))
\end{align*}
We then set up $\displaystyle g:\left(\mathbb{R}^{n}\right)^{2}\rightarrow \mathbb{R}^{2}$ as $\displaystyle g(\mathbf{x} ,\mathbf{y}) \mapsto B(( 0,0) ,(\mathbf{x} ,\mathbf{y}))$. This gives $\displaystyle g(\mathbf{x} ,\mathbf{y}) =( \Sigma _{i} a_{1}( x_{i}) \cdot y_{i} ,\Sigma _{i} a_{2}( x_{i}))$. We substitute $\displaystyle o=z=0$ into the above expression. This gives us;
\begin{align*}
B( g(\mathbf{x} ,\mathbf{y}) ,\ (\mathbf{x} ',\mathbf{y} ')) & =g(\mathbf{x} \oplus \mathbf{x} ',\mathbf{y} \oplus \mathbf{y} ')
\end{align*}
Showing that $\displaystyle g$ is streamable, interpreting the input as $\displaystyle \left(\mathbb{R}^{2}\right)^{n}$. We set up the tail $\displaystyle t:\mathbb{R}^{2}\rightarrow \mathbb{R}$ as $\displaystyle t( o,z) \mapsto o/b( z)$. Therefore, $\displaystyle t( g(\mathbf{x} ,\mathbf{y})) =( \Sigma _{i} a_{1}( x_{i}) \cdot y_{i}) /b( \Sigma _{i} a_{2}( x_{i})) =NC(\mathbf{x} ,\mathbf{y})$ is streamable.
\qed

\begin{figure}[h]
    \centering
\includegraphics[width=1\textwidth]{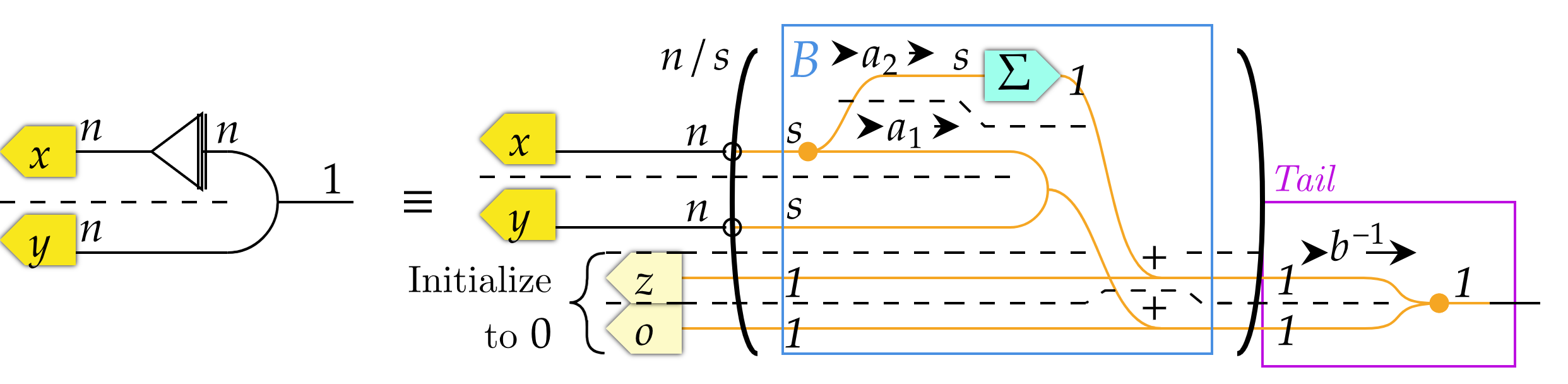}
    \caption{
    Normalized-contraction is mathematically expressed on the left-hand side, using a triangle with a triple line to indicate a generic normalization, followed by a dot product. The diagram reflects Lemma \ref{lem:normalized_contraction}, which deduces that a normalized contraction can be performed by a loop.
    }
    \label{fig:normalized_contractionproof}
\end{figure}

This lemma can be represented diagrammatically in Figure \ref{fig:normalized_contractionproof} by asserting that a normalized-contraction (left-hand side) returns the same result for the same inputs ($\equiv$) as a decomposed looped operation (right-hand side), which partitions data and runs the accumulator. Note that $\equiv$ does not assert that the computational cost of both sides are equivalent, merely that they preserve ``correctness''. As the computational costs may differ, $\equiv$ show how algorithms can be optimized to have desirable resource usage profiles.

\begin{theorem}
Attention where the SoftMax is replaced by another normalization operation remains streamable, as shown in the diagram below;
\end{theorem}

\begin{figure}[h]
    \centering
\includegraphics[width=1\textwidth]{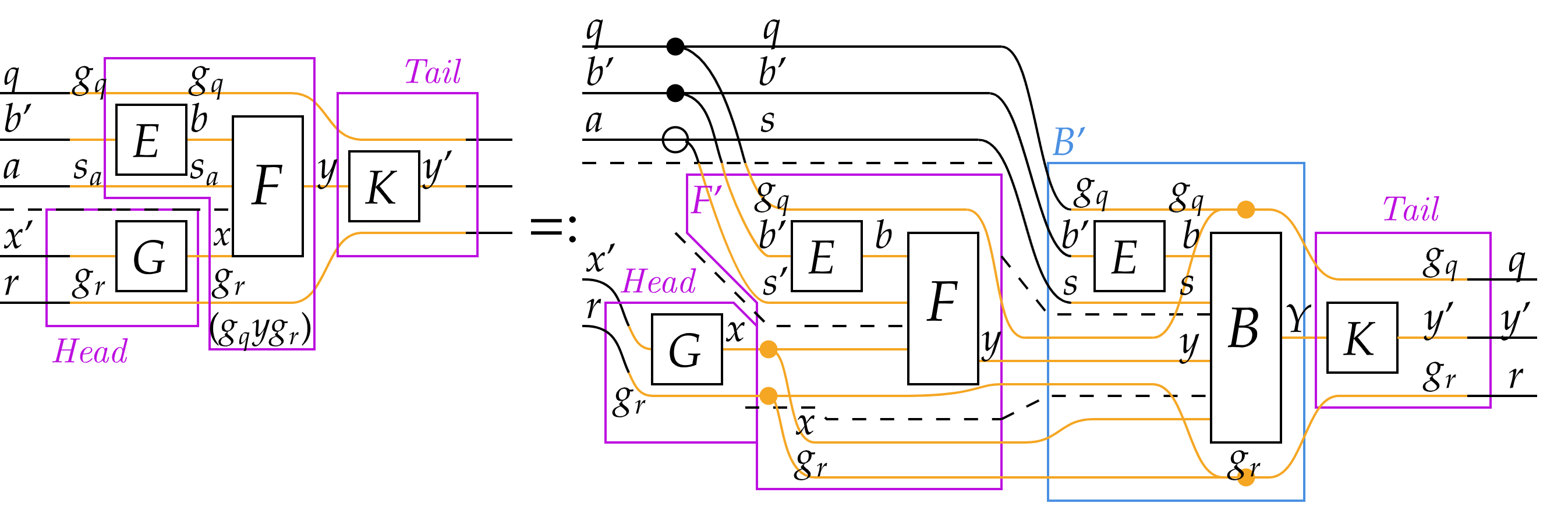}
    \caption{
    The fusion theorems indicate how modifying streamable operations changes the underlying accumulator, while maintaining streamability.
    }
    \label{fig:fusion_theorems}
\end{figure}

\textit{Proof.} From \cite{abbott_flashattention_2024}, we have the fusion theorem, which indicates how composed and broadcasted streamable operations maintain their streamability. It can be diagrammatically represented in Figure \ref{fig:fusion_theorems}.

We start with a diagram of streamable normalized-contraction. Then, we apply the streaming theorems to add the $\displaystyle Q/K$ contraction and broadcasting over the query and key/value axes. These operations preserve the streamable axis, and therefore provide a streamable algorithm. Figure \ref{fig:fusion_theorems} indicates how the new algorithm can be fully expressed with a loop.

\begin{figure}[h]
    \centering
\includegraphics[width=1\textwidth]{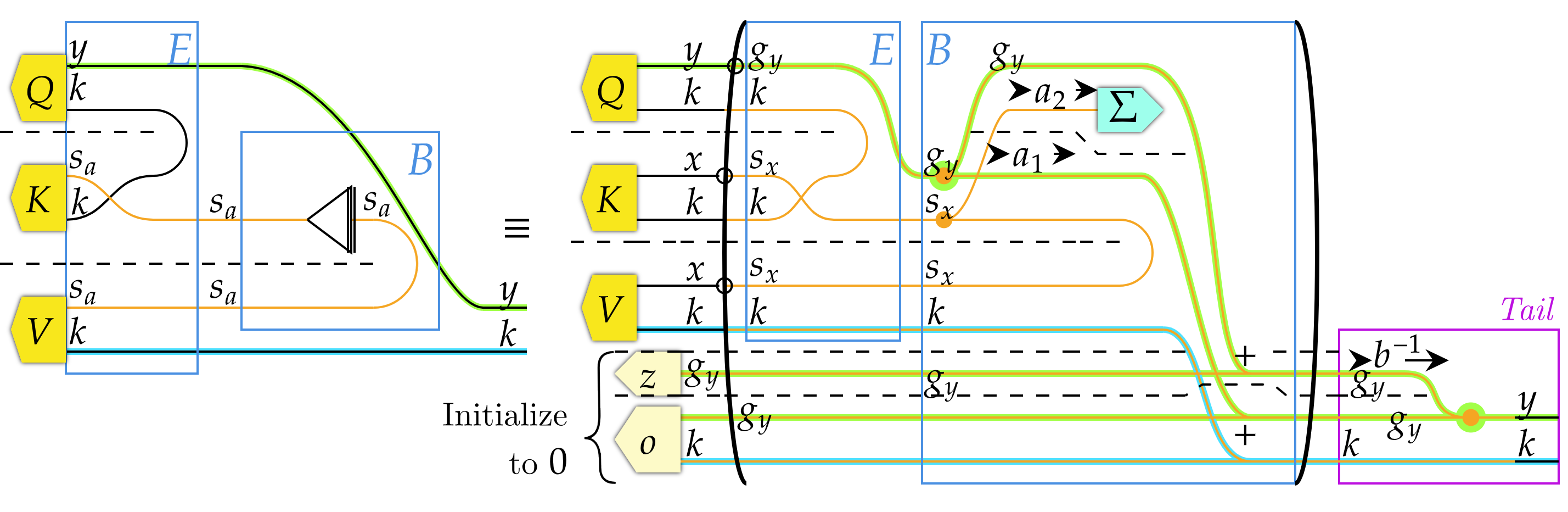}
    \caption{
    An attention algorithm using a generic normalization matches the left-hand side of the diagram. This matches Figure \ref{fig:normalized_contractionproof}, but with modifications similar to the fusion theorems, Figure \ref{fig:fusion_theorems}. This results in the fully expanded loop form on the right-hand side.
    }
    \label{fig:attention_expansion}
\end{figure}
\qed
\FloatBarrier

\subsection{Model Details}\label{apx:full_model}
The full model utilizes additional features, which match the standards of the latest models. RMSNorm (circle with a V-sign, to represent morning with variance) is the norm used in 
-V3 \cite{deepseek-ai_deepseek-v3_2025} and Mixtral-8x7B \cite{jiang_mixtral_2024}. Residual connections \cite{he_deep_2015,he_identity_2016} have been present since the original transformer architecture \cite{vaswani_attention_2017}. Grouped query attention is used \cite{ainslie_gqa_2023} as Multi-head Latent Attention (MLA) is not directly applicable, as we have no rotary, positional information due to the bag structure.

\subsection{Limitations}
\textbf{SMEM Memory Banking} SMEM data is stored in banks of 32 4B registers. SMEM/RMEM transfers from the same warp can be done in parallel, 128B transactions of variables in 32 distinct banks. However, we are required to load $\displaystyle Q$, $\displaystyle K$, and $\displaystyle V$ in tiles of size $\displaystyle 8\times 8$ to tensor cores. At $\displaystyle k=128$, an $\displaystyle 8\times 8$ tile straddles only $\displaystyle 8$ banks, reducing SMEM throughput by 75\%. By padding the end of each row so that $\displaystyle k'=136$, each $\displaystyle 8\times 8$ tile can be made to straddle all the banks. $\displaystyle k'=128+8$ was chosen as it maintains the 16-byte alignment required by the A100's cp.async.cg.shared.global PTX operation, which allows the L1 cache to be skipped. As diagrams indicate the cache requirements, we can safely skip it. The memory banking optimization is not currently represented by diagrams. Diagrams, however, already have information regarding axis size and stride. If information regarding the underlying tile and bank sizes can be incorporated, this optimization may be deduced.

\textbf{Numerical Accuracy} Though we have left an analysis of the numerical accuracy of this algorithm to future work, we tested against a reference FP32 implementation in PyTorch to ensure close results, with \(99.7\%\) of elements falling within 0.01 of the reference.

\textbf{Warp Shuffling} On the A100, tensor cores operate on fragmented memory. FP16 tensor core data is on threads, however, it is distributed in $8\times 8$ matrix fragments distributed across warps of 32 threads. This makes axis-wise operations, such as row accumulations, difficult to manage. Saving data to SMEM allows it to be coalesced. However, this creates pressure on the SMEM/RMEM transaction pipeline. Warp shuffles allow threads to directly exchange data, avoiding SMEM. This can improve performance, and is diagrammatically akin to blue and green levels directly transferring data. However, its implementation is complex, often requiring external packages such as \href{https://github.com/NVIDIA/cutlass}{CUTLASS} to be properly managed. Automating the development of a kernel from diagrams, however, will allow us to overcome this challenge as the complexities of how data is moved between blue and green levels is managed in the background.

\subsection{Experiment Details}
We ran our benchmarking loop for 100 warmup runs and 100 timed runs at each sequence length. We benchmarked our results on DataCrunch 1A100.22V machine consisting of 22 cores allocated from an AMD Epyc 7642 48-core CPU, 120GB RAM and 1 A100 SXM4 80GB GPU. 

\end{document}